\documentclass[12pt]{amsart}

\usepackage{graphicx}

\begin{document}


\title{Ribbonlength of Torus Knots}       
\author{Brooke Brennan, Thomas W. Mattman, Roberto
Raya, and Dan Tating}        
\address{Department of Mathematics and Statistics,
         California State University, Chico,
         Chico CA 95929-0525, USA}
\email{bbrennan@rbuhsd.k12.ca.us}
\email{tmattman@csuchico.edu}
\email{rraya@mail.csuchico.edu}
\email{dtating@tco.net}
\thanks{The first author is a high school teacher partially supported by 
NSF REU Award \# 0354174: Research Experiences in Mathematics for
Undergraduates and Teachers. The third and fourth
authors are undergraduate students who received support from the MAA's program for
Strengthening Underrepresented Minority Mathematics Achievement (SUMMA) with
funding from the NSF and NSA. These three were supervised by the second author who
also received support through the SUMMA and NSF REU programs. The 
research took place in part during an REU at CSU, Chico
during the summer of 2004.}

\begin{abstract}
Using Kauffman's model of flat knotted ribbons, we demonstrate how 
all regular polygons of at least seven sides can be 
realised by ribbon constructions of torus knots. 
We calculate length
to width ratios for these constructions thereby bounding the 
Ribbonlength of the knots. In particular, we give evidence that the
closed (respectively, truncation) Ribbonlength of a
$(q+1,q)$ torus knot is $(2q+1) \cot(\pi/(2q+1))$ (resp., $2q \cot(\pi/(2q+1))$).
Using these calculations, we provide the bounds $c_1 \leq 2/\pi$ and 
$c_2 \geq 5/3 \cot \pi/5$ for the constants $c_1$ and $c_2$
that relate Ribbonlength $R(K)$ and crossing number $C(K)$ in a conjecture of
Kusner: $c_1 C(K)  \leq R(K) \leq c_2 C(K)$. 
\end{abstract}

\maketitle

\markboth{BRENNAN, MATTMAN, RAYA, AND TATING}
{RIBBONLENGTH OF TORUS KNOTS}

\section{Introduction}

In \cite{K}, Kauffman introduces a model for knots presented as 
flat knotted ribbons and gives constructions of the trefoil and figure eight
knots. He defines the (closed) {\em Ribbonlength} of a knot to be the smallest
length to width ratio possible among the ways of forming the knot as a closed loop
of ribbon. A {\em truncation presentation} of a knot is one formed from a length of
ribbon such  that the ends of the ribbon lie flush with segments in its edges.
The minimum length to width ratio over such truncation presentations is the {\em
truncation Ribbonlength}.

In her Master's thesis,
DeMaranville~\cite{D} showed how to build regular
polygons by tying torus knots with ribbon. In the current paper, we summarise
some of the constructions of
\cite{D} and use them to bound the Ribbonlength of certain torus
knots. In particular, we conjecture that the closed (respectively,
truncation) Ribbonlength of a
$(q+1,q)$ torus knot is $(2q+1) \cot(\pi/(2q+1))$ (resp., $2q \cot(\pi/(2q+1))$).
We illustrate how to construct regular polygons by forming
$(q+1,q)$, $(p,2)$, $(2q+2,q)$, and
$(2q+4,q)$ torus knots with ribbons. In fact, these families
include all regular polygons except triangles, squares, and hexagons. 

Kauffman~\cite{K} reports that Kusner conjectures a linear relationship between
Ribbonlength and crossing number:
\begin{equation}
c_1 \mbox{Crossing} (K) \leq \mbox{Ribbonlength}(K) \leq c_2 \mbox{Crossing}
(K) 
\label{eqKusner}
\end{equation}
Using estimates of the Ribbonlength for torus knots, we deduce bounds on the 
constants $c_1$ and $c_2$. For truncated knots, we conclude
$c_1 \leq 4/ \pi$ while for closed knots we have $c_1 \leq 2 / \pi$.
For $c_2$, we cannot improve on the bounds that follow from Kauffman's~\cite{K}
study. The closed trefoil shows $c_2 \geq  \frac53 \cot(\pi / 5)$ while
the truncated figure eight knot yields $c_2 \geq (3 + \sqrt{2})/2$.

The paper is structured as follows. After this introduction we devote one
section to each of the following families of torus knots: $(q+1,q)$, 
$(p,2)$, $(2q+1,q)$, and $(2q+2,q)$ and $(2q+4,q)$. In each case we
determine the length to width ratio of the construction. Together, these families
yield all regular $n$-gons for $n > 6$ as well as the regular pentagon.
Next, in Section 6, we 
refine our estimates of the Ribbonlength for the $(5,2)$ and $(7,2)$ torus knots
and show how the $7_4$ knot (which is not torus) can be tied as a rectangle. 
In the final section we compile the length to width ratios of
the various families of torus knots in order to obtain bounds on the constants
$c_1$ and $c_2$ of Equation~\ref{eqKusner}.

\section{$(q,q+1)$ torus knots}

In \cite{K}, Kauffman uses the pentagon formed when a  
trefoil knot (a $(3,2)$ torus knot) is tied with a ribbon to find  the
(conjectured) Ribbonlength of that knot. DeMaranville~\cite{D} observed that a
$(4,3)$ torus knot forms a heptagon and, in general, a $(q+1,q)$ torus knot
yields a $2q+1$-gon. We will use this observation to calculate the Ribbonlength of
these knots.

\begin{figure}[ht]
\begin{center}
\includegraphics[scale=0.3]{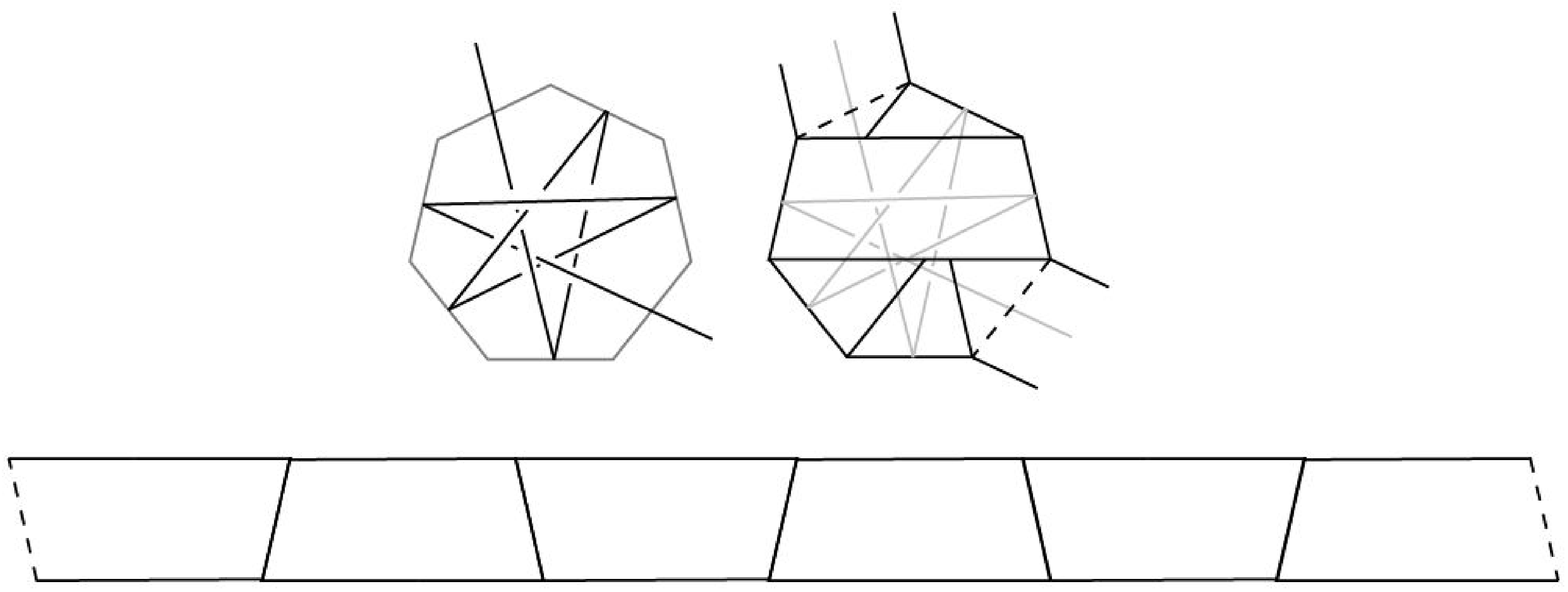}
\caption{The $(4,3)$ torus knot forms a heptagon. \label{fig43}}
\end{center}
\end{figure}

\begin{figure}[ht]
\begin{center}
\includegraphics[scale=0.3]{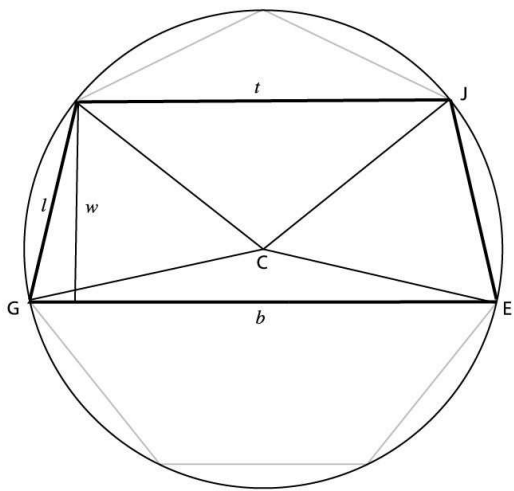}
\caption{The circumscribed heptagon. \label{fig43circ}}
\end{center}
\end{figure}

To illustrate the geometry of the $(q,q+1)$ torus knots, we begin with the
truncated $(4,3)$ torus knot. As in Figure~\ref{fig43},
when the ribbon is folded
to  make a regular heptagon, the fold lines segment the ribbon into six congruent
isosceles
trapezoids. By circumscribing the heptagon, we recognise the top and base of each
trapezoid as chords of a circle (Figure~\ref{fig43circ}).
We will use elementary geometry
to find those lengths.

In Figure~\ref{fig43circ}, $C$ is the centre of the circle.
Since $ \angle GCE = \pi - \pi/7$, the remaining angles in the
isosceles triangle
$\Delta GCE$ are each $\frac12 (\pi /7)$. By the Law of Sines, the length $b$ of
the base of the trapezoid is related to the radius $r$ by
$b = r \sin (\angle GCE) / \sin(\angle CEG)$. Since $\Delta CEJ$ is isosceles,
we can find $r$ in terms of $l$, the length of the sides of the heptagon:
$r = l/( 2 \sin(\frac12 (2 \pi/ 7)))$. Thus, 
\begin{eqnarray*}
b & = & \frac{l \sin (\pi - \pi/7 )}{2 \sin( \pi/7) \sin( \frac12 \pi/ 7)} \\
  & = & \frac{l}{2 \sin (\frac12 \pi/7)}
\end{eqnarray*}

Similarly, the top of the trapezoid has length
\begin{eqnarray*}
t & = & \frac{l \sin (\pi - 3 \pi/7)}{2 \sin (\pi /7) \sin ( \frac 32 \pi / 7)} \\
  & = & \frac{ 2 l \sin (\frac32 \pi /7) \cos(\frac32 \pi /7)}{2 \sin (\pi /7)
\sin( \frac32 \pi / 7)} \\
   & = & \frac{l \cos (\frac32 \pi /7)} {\sin (\pi /7)}
\end{eqnarray*}

The length of the ribbon is $L = 3(b+t)$. In order to compare it with the 
width $w$, note that $w = l \sin(3 \pi/7) = l \sin (\pi/2 - \frac12 \pi/7) = 
l \cos(\frac12 \pi /7)$. So, the length to width ratio of the $(4,3)$ torus knot
is 
\begin{eqnarray*}
L/w & = & \frac{3 l}{l \cos(\frac12 \pi/7)} \left[ \frac{1}{2 \sin(\frac12 \pi/7)}
+ \frac{ \cos(\frac32 \pi /7)}{\sin( \pi/7)} \right] \\
& = & \frac{3}{\cos(\frac12 \pi/7)} \left[ \frac{\cos (\frac12 \pi / 7)}{2
\sin(\frac12
\pi/7)\cos(\frac12 \pi / 7)} +
\frac{ \cos(\frac32 \pi /7)}{\sin( \pi/7)} \right] \\
& = &  \frac{3}{\cos(\frac12 \pi/7)} \left[ \frac{\cos (\frac12 \pi / 7)}{\sin(
\pi/7)} + \frac{ \cos(\frac32 \pi /7)}{\sin( \pi/7)} \right] \\
& = & \frac{3(2 \cos((\frac12 \pi/7+ \frac32
\pi/7)/2) \cos((\frac12 \pi/7 - \frac32 \pi/7)/2))}{\cos(\frac12 \pi/7)
\sin(\pi/7)}  \\
& = & \frac{6 \cos(\pi/7)
\cos(\frac12 \pi/7 )}{\cos(\frac12 \pi/7) \sin(\pi/7)} \\
& = &  6 \cot(\pi/7)
\end{eqnarray*}

We conjecture that $6 \cot(\pi/7)$ is the Ribbonlength of this knot; we expect
that there is no way to tie this knot with a ribbon that will result in a 
smaller length to width ratio.

DeMaranville observed that the pattern of pentagon formed by a $(3,2)$ torus
and heptagon by a $(4,3)$ torus persists. In general, a $(q+1,q)$ torus knot
can be tied with ribbon to form a regular $(2q+1)$-gon. The
calculation above also generalises. The ribbon will be segmented into $2q$
isosceles trapezoids with tops and bottoms of length $t = \frac{l \cos (\frac32 \pi
/(2q+1))} {\sin (\pi /(2q+1))}$ and $b = \frac{l}{2 \sin (\frac12 \pi/(2q+1))}$
respectively. Comparing the length $q(b+t)$ with the width $w = l \cos(\frac12
\pi/(2q+1))$ we arrive at our conjectured truncation Ribbonlength
for a $(q+1,q)$ torus knot: $2q \cot(\pi/(2q+1))$. For the closed knot, we
need one additional trapezoid. Measuring the length along the centre of the
ribbon, we conjecture that the closed Ribbonlength of a $(q+1,q)$ torus
knot is $(2q+1) \cot(\pi/(2q+1))$.

\section{$(p,2)$ torus knots}

\begin{figure}[ht]
\begin{center}
\includegraphics[scale=0.3]{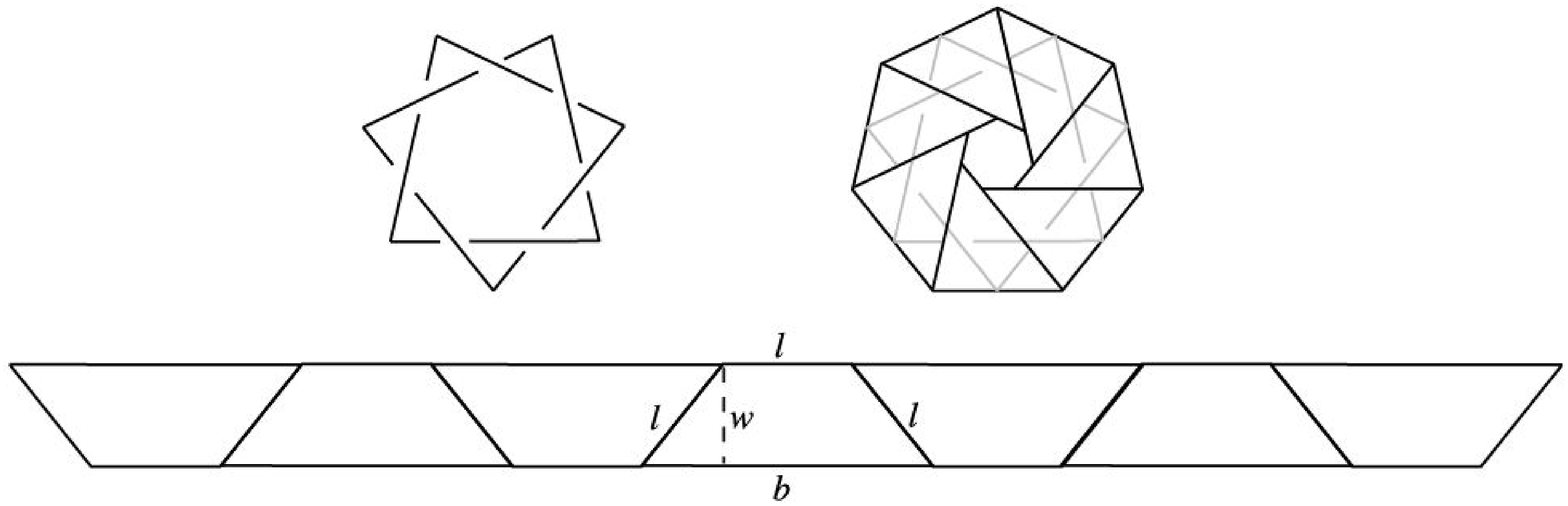}
\caption{The $(7,2)$ torus knot. \label{figp2}}
\end{center}
\end{figure}

In this section we 
investigate
$(p,2)$ torus knots ($p \geq 7$, odd). DeMaranville~\cite{D} proved that
knots in this family also form regular polygons when tied with ribbon
(see Figure~\ref{figp2}).
Note that there is a ``hole" in the centre of the
polygon that is in the shape of a smaller concentric regular $p$-gon.

As illustrated in Figure~\ref{figp2}, in tying this knot, the ribbon is segmented
into $p$ congruent isosceles trapezoids. Three of the sides of the trapezoid are
equal  to the side $l$ of the regular $p$-gon. The base of the trapezoid is again
a chord of the circumscribing circle of length
$b = l(1 + 2 \cos(2 \pi/ p))$. The length of the ribbon (measured along its
centre line) is $L = \frac{p}{2}(l+ b)$ while the width
is $w = l \sin(2 \pi/p)$. The closed length to width ratio for these knots is
therefore
\begin{eqnarray*}
L/w & = & \frac{2 p l (1 + \cos(2 \pi/p))}{2 l \sin(2 \pi/p)}  \\
   & = & \frac{p (1+ \cos^2(\pi/p) - \sin^2(\pi/p))}{2 \sin( \pi/p) \cos( \pi/p)}
\\
   & = & p \frac{2 \cos^2(\pi/p)}{2 \sin( \pi/p) \cos( \pi/p)} \\
   & = & p \cot(\pi/p).
\end{eqnarray*}

We anticipate that this length to width ratio is not the Ribbonlength for these
knots. For example, in Section~\ref{sec52} we discuss another method for tying
the $(7,2)$ torus knot that leads to a smaller length to width ratio than
the value $7 \cot(\pi/7)$ just obtained. On the other hand, the 
ratio of $p \cot(\pi/p)$ is an upperbound for the Ribbonlength of knots in
this family.

\section{$(2q+1,q)$ torus knots}

\begin{figure}[ht]
\begin{center}
\includegraphics[scale=0.27]{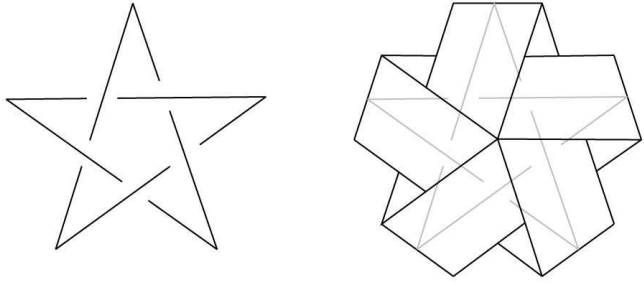}
\caption{The $(5,2)$ torus knot. \label{fig2qp1}}
\end{center}
\end{figure}

Unlike the previous two families, the $(2q+1,q)$ knots ($q > 1$) form ``pinwheels"
rather than polygons when tied with ribbon (Figure~\ref{fig2qp1}).
However, like the other
families,  the ribbon is segmented into congruent isosceles trapezoids. 
Each trapezoid has a top and base of length 
$t = 2w \cot(\pi/(2q+1))$ and 
$b = 2w(\cot(\pi/(2q+1)) +
\tan(\pi/2(2q+1))) = 2w/ \sin(\pi/(2q+1))$ respectively, where $w$ is the width of
the ribbon. The length to
width ratio is, therefore, 
\begin{eqnarray*}
L/w & = & \frac{2q+1}{2} (b+t)/w \\
   & = & (2q+1) (\cot(\pi/(2q+1)) + \csc(\pi/(2q+1)))\\
   & = & (2q+1) \cot(\pi/(2(2q+1)))
\end{eqnarray*}
Again, while this provides an upper bound for the Ribbonlength, we anticipate that
the Ribbonlength for these knots is not realised by this construction.
In particular, we demonstrate below (Section~\ref{sec52}) how the $(5,2)$ torus
knot can be tied with a smaller length to width ratio than the value $5
\cot(\pi/10)$ just obtained.

\section{Torus knots that form even-sided polygons}

\begin{figure}[ht]
\begin{center}
\includegraphics[scale=0.38]{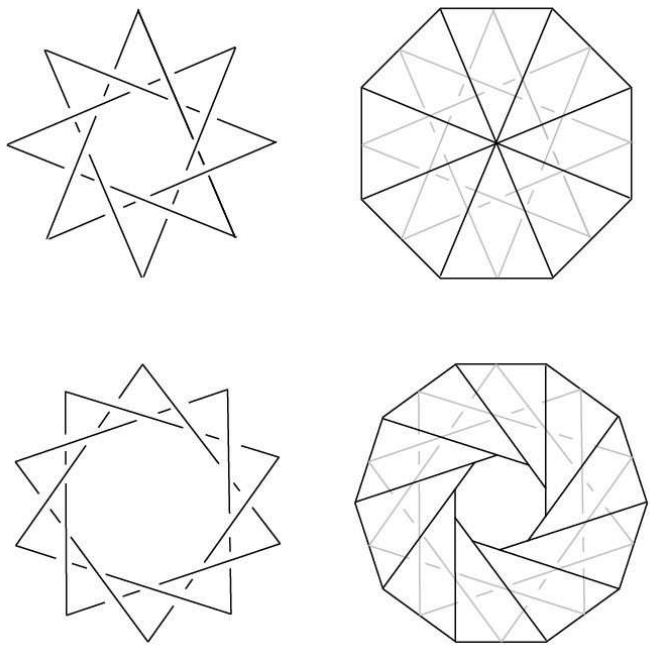}
\caption{The $(8,3)$ and $(10,3)$ torus knots. \label{figeven}}
\end{center}
\end{figure}

The regular polygons formed by $(q+1,q)$ and $(p,2)$ torus knots are odd-sided.
In this section we examine the $(2q+2,q)$ and $(2q+4,q)$ (where $q > 1$, odd) torus
knots which result in  $2q+2$- and $2q+4$-sided regular polygons~\cite{D},  see
Figure~\ref{figeven}. 
With the exception of the $(8,3)$ torus knot, these are
similar to the $(p,2)$ examples above in that there is a hole in centre of the
polygon.  They are also similar in that the length to width ratio is $n
\cot(\pi/n)$ where
$n = 2q+2$ or $2q+4$ is the number of sides of the polygon. Due to the 
polygonal hole in the centre, we anticipate that this length to width ratio is not
the Ribbonlength. The only exception is the $(8,3)$ torus knot which is formed
with no gap in the centre  (Figure~\ref{figeven}).
We conjecture that the Ribbonlength of that knot is $8 \cot(\pi/8)$.

The $(q+1,q)$ (for $q \geq 3$) and $(p,2)$ (for $p \geq 7$) torus knots give us two
different ways to construct each odd-sided polygon of seven or more sides. In
addition, the $(3,2)$ torus knot results in a regular pentagon. On the other hand,
the $(2q+2,q)$ and $(2q+4,q)$ knots (for $q \geq 3$)  yield all the even sided
regular polygons of eight or more sides. Together these families account for all
regular
$n$-gons of seven or more sides as well as the regular pentagon.

\section{The knot $7_4$ and the $(5,2)$ and $(7,2)$ torus knots
\label{sec52}}

In this section, we demonstrate ``shorter" ways of tying the $(5,2)$ 
and $(7,2)$ torus knots. This suggests that the length to width ratios
for the $(p,2)$ and $(2q+1,q)$ knots above are likely not the Ribbonlengths
of those knots. In addition, we look at how the $7_4$ knot results in a 
rectangle.

\subsection{The $(5,2)$ torus knot}

\begin{figure}[ht]
\begin{center}
\includegraphics[scale=0.27]{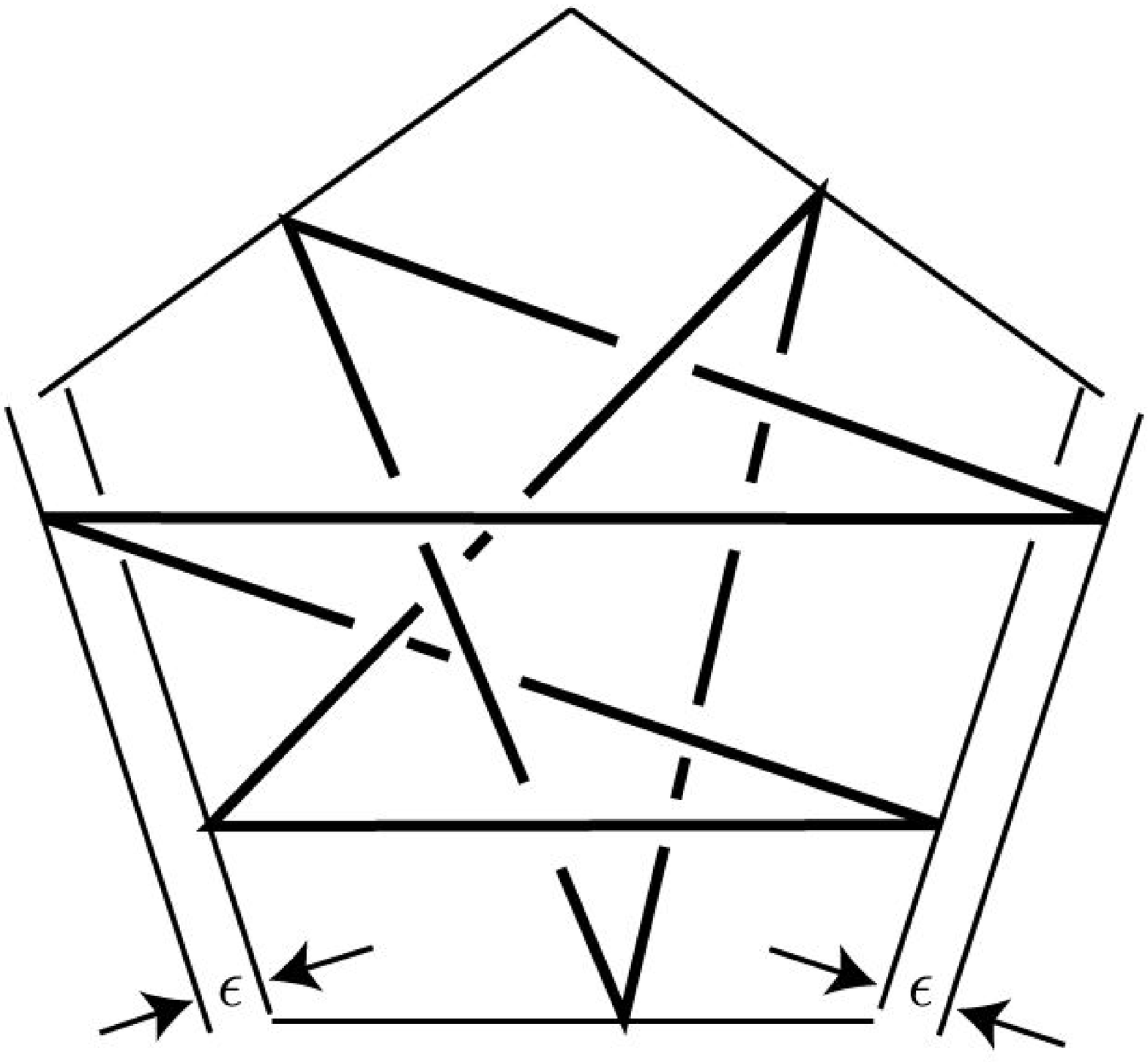}
\end{center}
\caption{A shorter version of the $(5,2)$ torus knot. \label{fig52}}
\end{figure}

We will tie a $(5,2)$ torus knot so as to obtain a smaller 
length to width ratio than the value $5 \cot( \pi/10) \approx 15.4$ obtained above.
In Figure~\ref{fig52},
we illustrate how to tie this knot to form a pentagon
much like the one given by the trefoil knot. Two sides of the pentagon are used
twice and we replace those sides by two parallel fold lines separated by a small
distance $\epsilon$. Although the three lines joining the doubled sides should
coincide, we have displaced them in the figure to show how those edges cross over
one another. 

Much like the trefoil knot, this way of folding $(5,2)$ will segment the ribbon
into $7$ trapezoids (as compared to the $5$ obtained for the closed trefoil). 
Therefore, this method of folding the $(5,2)$ torus will result in a length to
width ratio of approximately $7 \cot(\pi /5)$ or $9.6$.

\subsection{The $(7,2)$ torus knot}

\begin{figure}[ht]
\begin{center}
\includegraphics[scale=0.27]{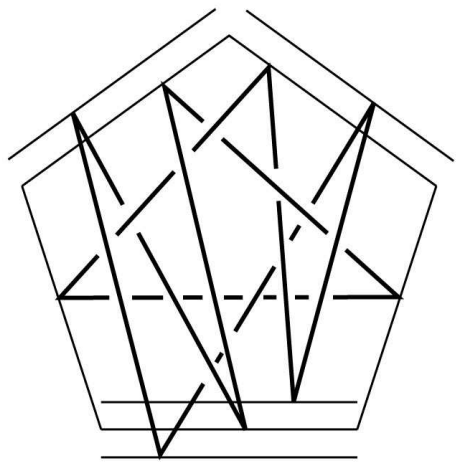}
\caption{A shorter version of the $(7,2)$ torus knot. \label{fig72}}
\end{center}
\end{figure}

The $(7,2)$ torus knot can also be tied to realise a length to width ratio
smaller than the value of $7 \cot(\pi/7) \approx 14.5$ mention above. As in
Figure~\ref{fig72}, 
this knot also results in a pentagon-like shape. 
(Again, coincident edges in the knot have been displaced in order to show
crossings.) This time
the ribbon is segmented into $9$ trapezoids resulting in a length to width ratio
of approximately $9 \cot( \pi/5)$ or $12.4$.

\subsection{The $7_4$ knot}

\begin{figure}[ht]
\begin{center}
\includegraphics[scale=0.27]{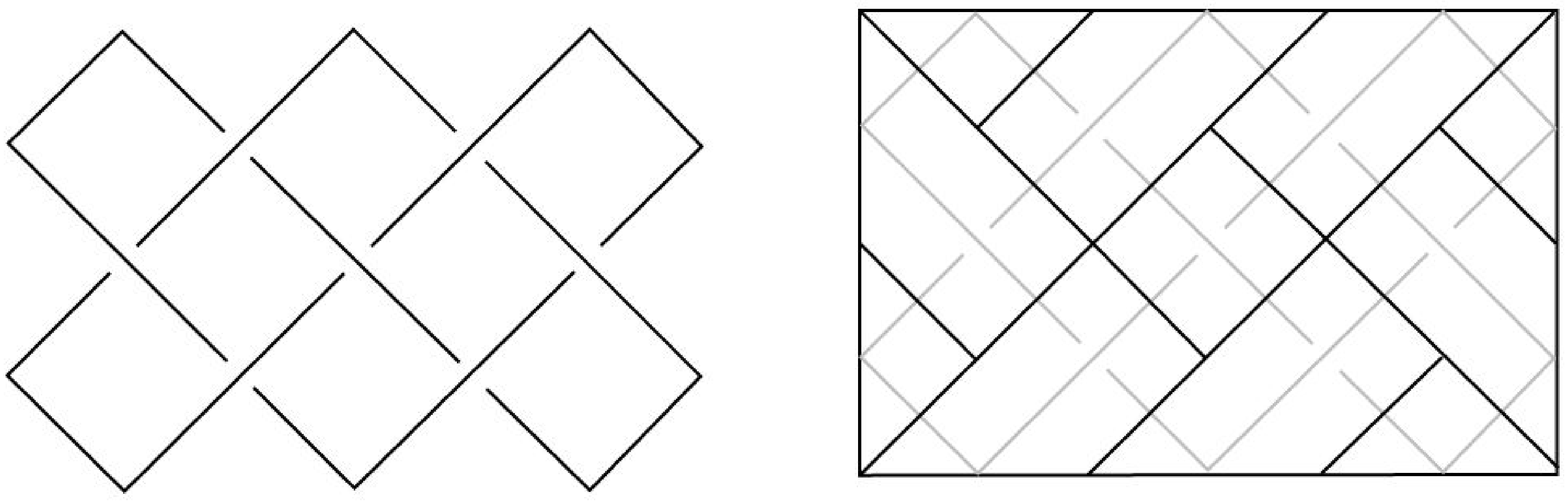}
\caption{The $7_4$ knot. \label{fig74}}
\end{center}
\end{figure}

Although it is not a torus knot, the $7_4$ knot can be folded to give a nice
geometric shape, a $2 \times 3$ rectangle~\cite{D}, see Figure~\ref{fig74}. 
The length to
width ratio of this pattern is $24$, as may easily be verified. However, this is 
likely not the Ribbonlength for this knot. Informal experiments suggest a
Ribbonlength of less than $20$.

\section{Ribbonlength and crossing number}

\setcounter{equation}{0}

In \cite{K}, Kauffman reports on Kusner's conjecture of a relationship between
Ribbonlength and crossing number:
\begin{equation}
c_1 \mbox{Crossing} (K) \leq \mbox{Ribbonlength}(K) \leq c_2 \mbox{Crossing}
(K)
\end{equation}
As the crossing number of a $(p,q)$ torus knot is $\min \{ p(q-1), q(p-1) \}$,
we can use our calculations of length to width ratios above to bound the 
constants $c_1$ and $c_2$. To this end, for each family of torus knots discussed
above, we divide the length to width ratio by the crossing number.

\subsection{$(q+1,q)$ torus knots}

Dividing the conjectured Ribbonlength by 
the crossing number $q^2-1$ yields 
$\frac{2q}{q^2-1} \cot(\pi/(2q+1))$ for the truncation presentation and
$\frac{2q+1}{q^2-1} \cot(\pi/(2q+1))$ for the closed knot. As $q$ increases, the
value of this ratio decreases, approaching $4/ \pi$ in the
limit. The largest values are obtained with the trefoil knot $(3,2)$, 
$\frac43 \cot(\pi/5) \approx 1.83$ and $\frac53 \cot(\pi / 5) \approx 2.29$
for the truncated and closed knot respectively.

\subsection{$(p,2)$ torus knots}

Dividing the length to width ratio by the crossing number $p$, we have 
$\cot(\pi/p)$. Notice that this tends to infinity as $p$ increases. We take this
as further evidence that this configuration does not realise the Ribbonlength for
these knots. Since $p \geq 7$, the smallest value we obtain is
$\cot(\pi/7) \approx 2.07$ for the closed $(7,2)$ torus knot.

\subsection{$(2q+1,q)$ torus knots}

Dividing the length to width ratio by the crossing number $(2q+1)(q-1)$, we have 
$\cot(\pi/(2(2q+1)))/(q-1)$ which again tends to $4/ \pi$. The largest value
obtained in this family is $\cot(\pi/10) \approx 3.1$ for the $(5,2)$ torus knot.

\subsection{Even-sided polygon knots, the figure eight knot, and $7_4$}

For the $(2q+2,q)$ torus knots, dividing the length to width ratio by
the crossing number $(2q+2)(q-1)$ yields $\cot(\pi/(2q+2))/(q-1)$ which tends
to $2/\pi$ as $q$ goes to infinity. The $(8,3)$ knot gives the biggest 
ratio in this family: $\cot(\pi/8)/2 \approx 1.2$. The crossing number of
the $(2q+4,q)$ knots is $(2q+4)(q-1)$. The resulting ratio
$\cot(\pi/(2q+4))/(q-1)$ again tends to $2 /\pi$ as $q$ goes to infinity. 
The largest value in this family is the $(10,3)$ knot for which
$\cot(\pi/10)/2 \approx 1.5$ 

For the $7_4$ knot,
on dividing the length to width ratio $24$ by the crossing number $7$, we obtain
$24/7 = 3 \frac37$. In \cite{K}, the Ribbonlength of the (truncated) figure eight
knot is conjectured to be $6+2 \sqrt{2}$. Dividing by the $4$ crossings, we have
$(3 +\sqrt{2})/2 \approx 2.2$.

\subsection{Bounding $c_1$ and $c_2$}

Each of the length to width ratios we have calculated above serves
as an upper bound for the Ribbonlength of the corresponding knot. 
Therefore, the comparisons of these ratios with crossing number
we have made in this section all serve to bound the constant $c_1$ of
Equation~\ref{eqKusner} above. For closed knots, we have the bound
$c_1 \leq 2/\pi$ realised by torus knots of the form $(2q+2,q)$ and
$(2q+4,q)$ in the limit as $q$ goes to infinity. For truncated knots,
the best bound is $c_1 \leq 4 / \pi$ realised by $(q+1,q)$ torus knots
as $q$ tends to infinity.

In order to estimate $c_2$, we restrict attention to those knots where we
believe we know the Ribbonlength, namely the $(q+1,q)$ torus knots, the knot
$(8,3)$ and the figure eight knot \cite{K}. For truncated knots the figure 
eight knot results in the  bound 
$c_2 \geq (3 + \sqrt{2})/2$. For closed knots, the trefoil yields the bound
$c_2 \geq \frac53 \cot(\pi / 5)$.


\end{document}